\documentclass[11pt]{article}
\usepackage{amsfonts,amsmath,fullpage, euscript}

\title {\bf  On Cover-Free Families
\thanks{Research supported by NSERC grant 239135-01}}
\author{R.\ Wei\\
  Department of Computer Science\\
 Lakehead University\\
 Thunder Bay, Ontario P7B 5E1\\
 Canada\\
  {\tt rwei@lakeheadu.ca}
 }

\date{April 22, 2006}
\newtheorem{Theorem}{Theorem}[section]
\newtheorem{Definition}{Definition}[section]

\newtheorem{Lemma}[Theorem]{Lemma}

\newtheorem{Corollary}[Theorem]{Corollary}
\def\whitebox{{\hbox{\hskip 1pt
        \vrule height 6pt depth 1.5pt
        \lower 1.5pt\vbox to 7.5pt{\hrule width
                  3.2pt\vfill\hrule width 3.2pt}%
        \vrule height 6pt depth 1.5pt
        \hskip 1pt } }}
\def\qed{\ifhmode\allowbreak\else\nobreak\fi\hfill\quad\nobreak
     \whitebox\medbreak}

\newcommand{\p}{{\noindent \it Proof.}\ \ }

\newcommand{\F}{{\cal F}}
\newcommand{\C}{{\cal C}}
\newcommand{\B}{{\cal B}}
\begin{document}
\maketitle
\begin{abstract}
Cover-free families were considered from different subjects by
numerous researchers. Recently, cover-free families are also found
useful in cryptography. In this paper, we use an uniform method
to survey known results of cover-free families. Many old results
are updated or generalized. Some new results are also given.
\end{abstract}

\section{Introduction}

Cover-free families were considered from different subjects such
as information theory, combinatorics and group testing. In 1964
cover-free families were first introduced   by Kautz and Singleton
\cite{ks} (they used a different terminology though) to
investigate nonrandom superimposed binary codes. These codes may
be used for retrieval files, data communication and magnetic
memories. After that, several researchers investigated in that
direction (see \cite{dr,drr,nz,kl,dswcp, dmr1,dmr}). In 1985,
Erd\"{o}s, Frankl and F\"{u}redi \cite{eff2} discussed cover-free
families as combinatorial objects which generalized the Sperner
systems. In 1987, Hwang and \'{S}os \cite{hs} defined cover-free
families for non-adaptive group testing. Since then, cover-free
families have been discussed in several equivalent formulations in
subjects such as information theory, combinatorics and group
testing by numerous researchers (see, for example,
\cite{bfpr,dr,e,eff,eff2,f,hs,kl,r,kbt,swz,hllt}). Recently,
cover-free families were also used for many cryptographic problems
such as frameproof codes and traceability schemes
(\cite{bs,cfn,cfn2,gsy,ssw,stw,sw,sw2}), broadcast encryption
(\cite{gsw,krs,sw3}, key storage (\cite{dfft,mp}), multi-receiver
authentication (\cite{snw}), etc.
  Since the researches of cover-free
families are from different disciplines,  different terminologies
were used through out the literatures and many results were
rediscovered several times. On the other hand, the definition of
cover-free families was generalized in several ways for different
applications. It is useful to have a comprehensive summary of
known results  about cover-free families so that repeated work can
be easily avoided in the future. Furthermore, many different
methods were used to investigate cover-free families. We may
compare these methods to find better ones or generalize these
methods for various cover-free families.

 The main purpose of this paper is to use a uniform method to investigate
known results of cover-free
families and their generalizations.
Basically, we will review and generalize the results of
existence and construction of cover-free families. Some new results of 
cover-free
families will  also be
given. We will not discuss applications of cover-free families
in this paper. However, this
survey paper will be useful for all the applications mentioned in the above
references.

In the rest of this section, we will give definitions of cover-free families
and equivalent objects of cover-free families. In section 2, we discuss
bounds of cover-free families. In section 3, we give various constructions of
cover-free families.

\subsection{Definitions and notations}

A {\em set system} is a pair $(X, \F)$, where $X$ is a set of points and $\F$ is
a set of subsets (called {\em blocks}) of $X$. A set system $(X,\F )$ is called
$k$-uniform if $|F|=k$ for each $F\in \F$. Throughout this paper, we will use
$N$
 to denote the cardinality of $X$ and use $T$ to
denote the cardinality of $\F$.

Cover-free families were defined in \cite{ks,eff2,hs} independently. We adapt
the following definition.

\begin{Definition}\label{d.r}
A set system $(X,\mathcal F)$ is called a {\em $r$-cover-free family} (or {\em
$r$-CFF})
provided that, for
any $r$ blocks $A_1,\cdots,A_r \in \mathcal{F}$ and any other block
$B_0  \in \mathcal{F}$,
we have
\[
B_0 \not\subseteq \bigcup_{j=1}^rA_j.
\]
\end{Definition}

The above definition has been  generalized in two different directions.
In 1988, Mitchell and Piper \cite{mp} defined the
concept of key distribution patterns,
which is in fact a generalized type of cover-free family. Their definition can
be described as follows.

\begin{Definition}\label{d.wr}
A set system $(X,\mathcal F)$ is called a {\em $(w,r)$-cover-free family} (or
{\em
$(w,r)$-CFF})
provided that, for
any $w$ blocks $B_1,\cdots,B_w \in \mathcal{F}$ and any other $r$ blocks
$A_1,\cdots,A_r \in \mathcal{F}$,
we have
\[
\bigcap_{i=1}^wB_i \not\subseteq \bigcup_{j=1}^rA_j.
\]
\end{Definition}

Note that the ``classical definition'' of cover-free family,
Definition \ref{d.r}, is the case $w=1$ of the above definition.

In 1989, Dyachkov {\it et al}
\cite{drr} considered another generalization of cover-free families when they
investigated superimposed distance codes. Their definition is as follows.

\begin{Definition}\label{d.rd}
A set system $(X,\mathcal F)$ is called a {\em $(r;d)$-cover-free family} (or
{\em
$(r;d)$-CFF})
provided that, for
any   block $B_0 \in \mathcal{F}$ and any other $r$ blocks
$A_1,\cdots,A_r \in \mathcal{F}$,
we have
\[
\left| B_0  \left\backslash\left( \bigcup_{j=1}^rA_j\right)\right| > d,\right.
\]
where $d$ is a positive integer.
\end{Definition}

Recently, Kumar {\it et al} \cite{krs} considered a set system $(X,\F )$, in
which
any   block $B_0 \in \mathcal{F}$ and any other $r$ blocks
$A_1,\cdots,A_r \in \mathcal{F}$
satisfy
\[
\left| B_0  \left\backslash\left( \bigcup_{j=1}^rA_j\right)\right| \ge
\alpha\left| B_0 \right|,\right.
\]
where $0<\alpha < 1$ is a fixed constant.

These structures are also used in \cite{gsw}.
Since they only considered the uniformed set system in
\cite{krs,gsw}, their
definition is in fact a special case of Definition \ref{d.rd}.

Now we give a more general definition of cover-free systems which covers all the
cases discussed so far.

\begin{Definition}\label{d.wrd}
A set system $(X,\mathcal F)$ is called a {\em $(w,r;d)$-cover-free family} (or
{\em
$(w,r;d)$-CFF})
provided that, for
any $w$ blocks $B_1,\cdots,B_w \in \mathcal{F}$ and any other $r$ blocks
$A_1,\cdots,A_r \in \mathcal{F}$,
we have
\[
\left|\left( \bigcap_{i=1}^wB_i \right)\left\backslash\left(
\bigcup_{j=1}^rA_j\right)\right| > d,\right.
\]
where $d$ is a non-negative integer.
\end{Definition}

Sometimes, we will use $(w,r;d)$-CFF$(N,T)$ to denote a cover-free family to
specify that $|X|=N$ and $|\F |=T$.

\subsection{Some equivalent objects}\label{sub.equi}

There are many equivalent objects of   cover-free families, such
as superimposed codes, key distribution patterns, non-adaptive
group testing algorithms, etc. In this subsection, we will
introduce two combinatorial objects, coverings of order-interval
hypergraph and disjunct set systems, which are also equivalent to
cover-free families. These two equivalent objects will be used in
later sections.

Let $[ n] =\{ 1,\cdots , n\}$ and let $P_{n;l,u} = \{ Y \subseteq
[n] : l \le |Y| \le u\}$, where $0< l < u < n$, ordered by
inclusion. Define a class of order-interval hypergraphs $G_{n;l,u}
= (P,E)$ as follows. Let the set of points $P = P_{n;l,u}$ and let
the set of edges $E$  be the maximal intervals, i.e.,
$$
E =  \{ I =\{ C \subseteq [n]:  Y_1\subseteq C \subseteq Y_2 \} :  |Y_1| = l,
|Y_2| =  u, Y_1 ,Y_2 \subseteq [n]  \}.
$$

Coverings of order-interval hypergraph were discussed in \cite{e}.

\begin{Definition}
A (point) covering of a hypergraph is a subset of points $S$ such that each edge
of the hypergraph contains at least one point of $S$.
\end{Definition}

 The
definition of a disjunct system was first defined in \cite{ssw}.


\begin{Definition}
\label{DS.defn}
{\rm
A set system $(X, \mathcal{B})$ is
an {\it $(i,j)$-disjunct system} provided that,
for any $P,Q \subseteq X$ such that
$|P| \le i$, $|Q| \le j$ and  $P\cap Q =
\emptyset$, there exists a $B \in \mathcal{B}$
such that $P \subseteq B$ and
$Q \cap B = \emptyset$.
An $(i,j)$-disjunct system, $(X, \mathcal{B})$,
will be denoted as an $(i,j)$-DS$(v,b)$ if $|X| = v$
and $|\mathcal{B}| = b$.
}
\end{Definition}

Disjunct systems for  $i = 1$ were discussed in \cite{dh,bfpr}.
Definition \ref{DS.defn} is a generalization. The equivalence of
these objects are proved in the following theorems.

\begin{Theorem}
There exists a covering of $G_{n,l,u}$ of size $b$ if and only if there exists a
$(l,n-u)$-DS$(n,b)$.
\end{Theorem}
\p
$S$ is a covering of $G_{n,l,u}$  if and only if for  any $Y_1, Y_2 \subseteq 
[n],
 Y_1\subset Y_2,
|Y_1|=l, |Y_2|= u$, there is a $C \in S$ such that $Y_1\subseteq C
\subseteq Y_2$. This is equivalent to that for any $Y_1, Y_3
\subseteq [n], |Y_1|=l, |Y_3|= n- u, Y_1\cap Y_3 =\emptyset$,
there is some $C \in S$ such that $Y_1\subseteq C$ and $Y_3\cap C
= \emptyset$. \qed

A set system can be described by an incidence matrix.
Let $(X, \mathcal{B})$ be a set system where
$X = \{x_1, x_2, \dots , x_v\}$ and
$\mathcal{B} = \{B_1, B_2, \dots , B_b\}$.
The {\em incidence matrix} of $(X, \mathcal{B})$
is the $b\times v$
matrix $A = (a_{ij})$, where
\[a_{ij}=\left\{\begin{array}{ll}
1 & \mbox{if $x_j \in B_i$}\\
0 & \mbox{if $x_j \not\in B_i$.}
\end{array} \right.
\]
Conversely, given an incidence matrix, we can define an
associated set system in an obvious way.

Disjunct systems and cover-free families are dual incidence
structures. If $A$ is an incidence matrix of a disjunct system,
the $A^T$, the transpose of $A$, is an incidence matrix of a
cover-free family. Thus we have the following theorem.

\begin{Theorem}
\label{equivCFF-DS}
There exists an $(i,j)$-CFF$(N,T)$
if and only if there
exists an $(i,j)$-DS$(T,N)$.
\end{Theorem}

\section{Bounds}

It is easy to see that there is a trade-off between the values of $N$ and $T$ in
a cover-free family. We want to maximize the value of $T$ or minimize the value
of $N$ in case that other parameters are given. The bounds of $N$ or $T$ in
cover-free families were discussed by many researchers. Some authors considered
the bound of the rate $\log T\over N$. In this section we discuss these bounds.

\subsection{NP-completeness}

The first question about the bound of CFF is whether we can find an efficient
algorithm to determine the bound. Unfortunately, this is not the case. We
introduce some result from \cite{be,be2}, which tells us that in general to find
a bound of CFF is  difficult.

Let $G_{n,l,u}$ be an order-interval hypergraph. Define
\[
\tau (G_{n,l,u}) = \min \{ |S| : S \mbox{ is a covering of } G_{n,l,u} \}.
\]
By the results of Subsection \ref{sub.equi}, $\tau$ is an upper bound of $N$ for
certain $(i,j)$-CFF.

The following theorem in \cite{be2} shows that to find a value of $\tau$ is a
difficult problem.

\begin{Theorem}\cite{be2}\label{t.NP}
The problem of deciding $\tau(G_{n,l,u}) \le k$ (with input $G_{n,l,u}$ and $k$)
is NP-complete.
\end{Theorem}

The main idea of the proof of Theorem \ref{t.NP} is to reduce the problem to an
edge cover problem, which was proved to be NP-complete in \cite{be}.

\subsection{Bounds for $r$-CFF}

Now we consider the bounds for CFF. We first give an upper bound of $T$
in an uniform CFF, which was proved in
\cite{eff2}. The following lemma will be used  to prove that bound. This lemma
is a dual version of a lemma proved in \cite{frankl}.

\begin{Lemma}\cite{frankl}\label{l.frankl}
Suppose $(X,{\cal B})$ is a $t$-uniform set system, where $|X|=k$. If for any
$r$ blocks $B_1,\cdots, B_r \in {\cal B}$, we have $|\cup_{i=1}^r B_i|<k$ and
$rt
\ge k$, then $|{\cal B}|\le {k-1 \choose t}$.
\end{Lemma}

Using the above lemma, we can prove the following bound.

\begin{Theorem}\cite{eff2}\label{t.uniform}
In a $k$-uniform $r$-CFF,
\[
T \le {N\choose {\lceil {k\over r}\rceil}}\left/ {k-1 \choose {\lceil {k\over
r}\rceil-1}}\right. .
\]
\end{Theorem}
\p
Suppose $(X,{\cal F})$ is a $k$-uniform $r$-CFF. For $F\in {\cal F}$, define
\[
{\cal N}_t(F)=\{ T\subset F : |T|=t, \exists F'\neq F, F'\in {\cal F},T \subset
F'\}
\]
where $t=\lceil {k\over r}\rceil$. Let $T_1, \cdots , T_r \in {\cal N}_t$. Then
we have $|\cup_{i=1}^r T_i|\le k-1$ by the property of ${\cal N}_t$. Since
$rt\ge k$,
$$|{\cal N}_t|\le {k-1 \choose t}$$
 by Lemma \ref{l.frankl}. Therefore for each $F \in {\cal F}$, there are at
least
\[
{ k\choose t}-{k-1\choose t} = {k-1\choose t-1}
\]
subsets $T\subset F$, which are not contained in any $F'\in \F \backslash \{ F
\}$. It follows that
\[
|\F | {k-1\choose t-1} \le {N\choose t}.
\]
The conclusion follows.
\qed

For some special cases, the bound in Theorem \ref{t.uniform} can be reached. For
example, using probabilistic method, \cite{eff2} constructed $2r$-uniform
$r$-CFF with
\[
T = (N^2/(4r-2) ) - o(N^2).
\]
However, in general we still do not know whether the bound in Theorem
\ref{t.uniform} is the best possible bound.

 Bounds for an $r$-CFF without block size restriction have been studied by
several
researchers. Sperner's theorem in \cite{sperner} states that in a $1$-CFF
\[
T \le   {N \choose \left\lfloor {N \over 2} \right\rfloor}
 \]
for all positive integers $N \geq 2$. This bound is tight.

For $r\ge 2$, the bound
\[
{\log T\over N} \le {c\over r}
\]
for some constant $c$ is proved in \cite{nz,hs,eff2}. There were many
improvements for this bound.
The best known lower bounds on $N$ for general $r$ are found in
\cite{dr,f,r},
where different proofs of the following theorem can be found.

\begin{Theorem}\label{w=1}
For any $r \geq 2$, it holds that
\[
N \ge c{r^2\over \log r}\log T,
\]
for some constant $c$.
\end{Theorem}

The constant $c$ in Theorem \ref{w=1} is shown to be
approximately $1/2$ in \cite{dr},
approximately $1/4$ in \cite{f} and
approximately $1/8$ in \cite{r}. The simplest proof is given in \cite{f}, which
is modified from the proof of Theorem \ref{t.uniform}. We will use a similar
method as used in \cite{f} to prove a bound for $(r;d)$-CFF later.

\subsection{Bounds for $(r;d)$-CFF}

 An $(r;0)$-CFF is an $r$-CFF. A bound of $(r;d)$-CFF for $d\ge 1$ is first
considered in \cite{drr}. In this subsection, we first state their bound.

Let
\[
h(\alpha ) = -\alpha\log \alpha - ( 1-\alpha)\log (1-\alpha)
\]
be binary entropy and
\[
e_r={r^r \over (r+1)^{r+1}}, \ \ r = 1,2,\cdots .
\]
Define $e = {d\over N}$. Then the following theorem is proved in \cite{drr}.

\begin{Theorem}\label{t.drrbound}\cite{drr}
In a $(r;d)$-CFF(N,T),
\[
{\log T\over N} \le U_r(e),
\]
where, for $e\ge e_r$, the function $U_r(e)=0$ and for $0<e<e_r$ the sequence
$U_r=U_r(e)$ is defined recurrently:
\begin{enumerate}
\item
if $r=1$, then
\[
U_1=\left\{
\begin{array}{ll}
h\left({1\over 2}(1-\sqrt{8e(1-2e)}\right), & \mbox{ if \ \   } 0<e<e_1={1\over
4}\\
0,& \mbox{ if \ \ } e\ge {1\over 4}
\end{array}
\right.
\]
\item
if $r\ge 2$, then
\[
U_r=\min \{ 1-e/e_r,U_1/r,V_r(e)\},
\]
where $V_r$ is the unique solution of the equation
\[
V_r=\max_{v} \left\{ h\left({v\over r}\right)-(v+e)h\left({v\over
(v+e)r}\right)\right\}.
\]
The above maximum is taken over all $v$ for which
\[
0\le v \le 1-{V_r\over U_{r-1}}-e.
\]
\end{enumerate}
\end{Theorem}

The bound displayed above is  complicated and the
proof in \cite{drr} is rather involved. In the following we give a simple bound
for general
$(r;d)$-CFF. The proof of the bound is modified from a proof used in \cite{f}.

\begin{Theorem}\label{t.gbound}
In a $(r;d)$-CFF,
\[
T< r+{N \choose m}
\]
where $m = \lceil {2(N-r-d(r+1))\over r(r+1)}\rceil$.
\end{Theorem}
\p Suppose $(X, \F)$ is a $(r;d)$-CFF, where $|X|=N$. For a
positive integer $t \le n/2$, define
\[
\F_t =\{F\in \F: \exists A \subseteq F \mbox{ such that } |A|=t, A\not \subseteq
F' \mbox{ for any } F'\in \F\backslash \{ F\}\}
\]
and ${\cal A}$ be the family of all these $t$-subsets. Let
\[
\F_0=\{ F\in \F : |F|<t\}, \ \ {\cal B}=\{ T\subset X: |T|=t \mbox{ and }
\exists F\in \F_0\mbox{ such that } F\subset T\}
\]
Then
\[
|\F_0\cup \F_t|\le |{\cal A}|+|{\cal B}| \le {N \choose t}.
\]
Define $\F'=\F \backslash (\F_t \cup \F_0)$. Then for any $F\in \F'$ and $F_1,
F_2, \cdots, F_i \in \F\backslash \{ F\}$, where $0\le i \le r$, we have
\[
\left| F\left\backslash \bigcup_{j=1}^iF_j\right| >t(r-i)+d\right. .
\]
In fact, if
\[
\left| F\left\backslash \bigcup_{j=1}^iF_j\right| \le t(r-i)+d\right. ,
\]
then we can write $F=D\cup F_1\cdots \cup F_i\cup A_{i+1} \cdots \cup A_r$,
where $|D| = d, |A_j|= t, i+1\le j \le r$. Since $F\in \F'$, there exits an
$F_j\in \F\backslash \{ F\}$ such that $A_j \subseteq F_j$. Therefore
$|F\backslash (\cup_{l=1}^r F_l)| \le |D| = d$, a contradiction.

Now for $F_0,F_1, \cdots, F_r \in \F'$, we have
$|\cup_{i=0}^r F_i|= |F_0|+|F_1\backslash F_0|+\cdots+|F_r\backslash (F_0\cup
\cdots \cup F_{r-1})| \ge t{r(r+1)\over 2}+(d+1)(r+1)$. Let
\[
t=\left\lceil {2((N-r)-d(r+1))\over r(r+1)}\right\rceil.
\]
Then $|\cup_{i=0}^r F_i|\ge N+1$, which is impossible. Thus $|\F' |\le r$ and
the conclusion follows.
\qed

When $d=0$, we obtain the following theorem of \cite{f}.
\begin{Theorem}\cite{f}
In an $r$-CFF,
\[
T< r+{N \choose m}
\]
where $m = \lceil {2(N-r)\over r(r+1)}\rceil$.
\end{Theorem}

In \cite{drr}, the authors indicated that one can get an asymptotical upper
bound of $\log T\over N$ for a
$(r;d)$-CFF from Theorem \ref{t.drrbound} as follows:
\[
\left(2-{dr\over N}\right){\log r\over r^2}.
\]

  From Theorem \ref{t.gbound} we can obtain
the following asymptotical bound, which is slightly bigger than the bound given
in \cite{drr}.
\[
4\left(1-{dr\over N}\right){\log r\over r^2}.
\]

In \cite{b}, better bounds of $(r;d)$-CFF
for  special cases when $ r=1$ or $2$ are
given. For example, the bound for a $(2;d)$-CFF is given as follows.

\begin{Theorem}\cite{b}
In a $(2;d)$-CFF,
\[
T< {N \choose t^*}\left/ {2t^*+d-1\choose t^*}\right.
\]
where $t^* $ is the least integer value of $t$ such that
\[
N\le 5t+2+{d(d-1)\over t+d}.
\]
\end{Theorem}

\subsection{Bounds for $(w,r)$-CFF}

Let $N((w,r),T)$ denote the minimum value of $N$ in a $(w,r)$-CFF.
Lower bounds on $N((w,r),T)$ for arbitrary $w$ were first discussed by Dyer {\it
et al} \cite{dfft} in 1995. They proved the following bound for $N$.

\begin{Theorem}\cite{dfft}\label{t.dfft}
\[
N((w,r),T) \ge r(w\log T -\log r -w\log w)
\]
\end{Theorem}

In 1996, Engel \cite{e} proved the following two improved bounds. The first
bound
is implied in his proofs.

\begin{Theorem}\label{engel1}\cite{e}
\[
N((w,r),T) \ge {w+r-1\choose w} \log (T-r-w+2).
\]
\end{Theorem}

The second bound is an asymptotically one.

\begin{Theorem}\label{engel}\cite{e}
For any $\epsilon > 0$, it holds that
\[
N((w,r), T) \ge (1-\epsilon ) {(w+r-2)^{w+r-2}\over(w-1)^{w-1}(r-1)^{r-1}}\log
(T-r-w+2)
\]
for all sufficiently large $T$.
\end{Theorem}

The proofs in \cite{e} are based on some fractional matchings and fractional
coverings of order-interval hypergraphs.

Note that a $(1,r)$-CFF is an $r$-CFF. When $w=1$, the lower bound in Theorem
\ref{t.dfft} is $r(\log T - \log r)$ and the lower bound in Theorem \ref{engel1}
is $(r-1)\log (T-r+1)$. These bounds are weaker than the bound given in Theorem
\ref{w=1}.

Recently, two bounds for $(w,r)$-CFF were given in \cite{swz}. The proofs
of these bounds are based on the following lemma.

\begin{Lemma}\label{reduceN}
\[
N((w,r),T)\ge N((w,r-1),T-1)+N((w-1,r),T-1).
\]
\end{Lemma}

The first bound of \cite{swz} is as follows.

\begin{Theorem}\label{Nbound2}\cite{swz}
For $w,r \geq 1$ and $T \geq w+r > 2$, we have
\[
N((w,r),T)
\ge 2c{{w+r\choose w}\over \log (w+r)}\log T,
\]
where $c$ is the same constant as in Theorem \ref{w=1}.
\end{Theorem}

The second bound of  \cite{swz} is also an asymptotical bound.

\begin{Theorem}\label{Nbound3}\cite{swz}
For any integers $w, r \geq 1$, there exists an integer $T_{w,r}$ such that
\[
N((w,r),T) \ge 0.7c{{w+r \choose w}(w+r)\over \log {w+r\choose w}}\log T
\]
for all $T > T_{w,r}$, where $c$ is the same constant as in Theorem
\ref{w=1}.
\end{Theorem}

The bounds in Theorems \ref{Nbound2} and \ref{Nbound3} are usually better than
the bounds given in \cite{dfft,e}.
In particular, these bounds build on the bound of Theorem \ref{w=1}.
If we fix $w$ (respectively, $r$) and let $r$ (respectively, $w$) vary,
then these bounds are always stronger than the previous bounds.
For fixed $r$, Theorem \ref{Nbound2} improves Theorem \ref{engel1} by a factor
of $w / \log w$, and Theorem \ref{Nbound3} improves Theorem \ref{engel} by a
factor
of $w^2 / \log w$.
On the other hand, when $w=r$, Theorem \ref{engel} is better than  Theorem
\ref{Nbound3}
by a factor of $\sqrt{w}$.

\subsection{Bound for $(w, r; d)$-CFF}\label{remark}

In \cite{drr}, the following bound for $(1,r;d)$-CFF was proven:

\begin{Theorem}\label{t.1rdbound}
If $r>1$ and $d \geq 1$ are integers, then
\[
N((1,r;d),T) \ge c\left({r^2\over \log r}\log T + (d-1)r\right),
\]
where $c$ is some constant.
\end{Theorem}

\cite{sw2} proved the following bounds.

\begin{Theorem}
Suppose $r,w$ and $d$ are inteegers, $r>w \geq 1$,
and $d \geq 1$. Then
\begin{eqnarray*}
N((w,r;d),T) &\ge &c \, 4^{w-1}\left( 1-{1\over T}\right)\left(1 -{1\over 
T-2}\right) \cdots \left( 1-{1\over T-2w+2}\right)\\
 &  & \times \left({(r-w+1)^2\over \log (r-w +1)}\log (T -2w) + (d-1)(r-w+1)\right)
\end{eqnarray*}
for some constant $c$.
\end{Theorem}

\begin{Theorem}\label{Nbound2}
For $w,r \geq 1$ and $T \geq w+r > 2$, it holds that 
\[
N((w,r;d),T)
\ge 2c{ \binom{w+r}{w}\over \log (w+r)}\log T +{1\over 2}c \binom{w+r}{w}(d-1),
\]
where $c$ is the same constant as in Theorem \ref{t.1rdbound}.
\end{Theorem}

  \begin{Theorem}\label{Nbound3}
For any integers $w, r \geq 1$, there exists an integer $T_{w,r}$ such that 
\[
N((w,r;d),T) \ge 0.7c{ \binom{w+r}{w}(w+r)\over \log \binom{w+r}{w}}\log T + 
{1\over 2}c \binom{w+r}{w}(d -1)
\]
for all $T > T_{w,r}$, where $c$ is the same constant as in Theorem 
\ref{t.1rdbound}.
\end{Theorem}

\section{Constructions}

For the purpose of applications, researchers tried to find
efficient constructions of cover-free families. Three basic
methods were used to construct CFF: combinatorial methods, methods
from coding theory, and probabilistic method. In this section, we
review these constructions.

\subsection{Constructions from combinatorial designs}

In \cite{bfpr,eff2,ks,stw,sw}, combinatorial designs were used to
construct $r$-CFF. For the details of combinatorial designs we
used here, readers are refereed to \cite{cd}. The following
construction of $r$-CFF using $t$-designs can be found in
\cite{eff2,ks,sw}. The construction displayed here is  generalized
to $(r;d)$-CFF from $r$-CFF. First we give the definition of a
$t$-packing design as follows.

\begin{Definition}
A  {\em $t$-$(v, k,\lambda)$ packing design} is a
set system $(X, \mathcal{B})$, where $|X|$ = $v$,
$|B| = k$ for every $B \in \mathcal{B}$,
and every $t$-subset of $X$ occurs in at most $\lambda $
blocks in $\mathcal{B}$.
\end{Definition}

A $t$-packing design is an $r$-CFF for certain value of $r$. We
obtain the following construction.

\begin{Theorem} \label{thm:packingFPC}
If there exists a $t$-$(v,k,1)$ packing design
having $b$ blocks, then there exists a
$(r;d)$-CFF$(v,b)$, where
$r = \lfloor (k-d-1)/(t-1)\rfloor $.
\end{Theorem}

In \cite{sw}, orthogonal arrays were used to obtain $t$-packing
designs. An {\em orthogonal array} OA$(t,k,s)$ is a $k\times s^t$
array, with entries from a set of $s \ge 2$ symbols, such that in
any $t$ rows, every  $t\times 1$ column vector appears exactly
once.

Suppose  $\{ s_1, s_2, \cdots,
s_k\}$ is a column in an OA$(t,k,s)$. Define a block as
$$\{ (s_1,1), (s_2,2), \cdots, (s_k,k)\}$$
accordingly. In this way, we can obtain a $t$-$(ks,k,1)$ packing
design from an OA$(t,k,s)$.

It is known that
if $q$ is a prime power and $t < q$, then there exists an
OA$(t,q+1,q)$ (see \cite{cd}), and hence a
$t$-$\left( q^2+q, q+1, 1\right)$ packing design with $q^t$ blocks
exists.

Therefore we have the following result.

\begin{Theorem} For any prime power $q$ and any integer $t < q$,
there exists a $(\left\lfloor {q-d \over t-1} \right\rfloor;d)$-CFF$(q^2+q,q^t)$
\end{Theorem}

A $t$-$(v,k,1)$ design is a special case of $t$-$(v,k,1)$ packing design, in
which every $t$-subset occurs in exactly one block. In a $t$-design, the number
of blocks is
\[
{v\choose t}\left/ {k\choose t}\right. .
\]
There are many known results on the existence and constructions of $t$-$(v,k,1)$
designs. However, no $t$-$(v,k,1)$ design with $v>k>t$ is known to exist for
$t\ge 6$ (See \cite{cd} for the details).

In \cite{stw}, another combinatorial method is used to construct $(w,r)$-CFF,
which uses separating hash families. We can use this method to construct
$(w,r;d)$-CFF.

\begin{Definition}
\label{SHF.defn}
{\rm An {\em $(n,m,\{w_1,w_2\})$-separating hash family}
is a set of functions
$\mathcal{F}$, such that $|Y|=n$, $|X|=m$,
$f: Y \rightarrow X$
for each $f \in \mathcal{F}$, and for any
$C_1,C_2 \subseteq \{1, 2, \dots , n\}$ such that
$|C_1| = w_1$, $|C_2| = w_2$ and $C_1 \cap C_2 = \emptyset$,
there exists at least one
$f \in \mathcal{F}$ such that
\[ \{f(y) : y \in C_1\} \cap \{f(y) : y \in C_2\} = \emptyset.\]
The notation
SHF$(N; n, m, \{w_1,w_2\})$
will be used to denote an $(n,m,\{w_1,w_2\})$-separating hash family
with $|\mathcal{F}| = N$.
}
\end{Definition}

An SHF$(N; n, m, \{w_1,w_2\})$
can be depicted as an $N \times n$ matrix with entries
from $\{1, 2, \dots , m\}$, such that in any two disjoint sets $C_1$ and
$C_2$
of $w_1$ and $w_2$ columns (respectively), there exists
at least one row such that the entries in the columns of $C_1$
are distinct from the entries in the columns of $C_2$.

Now suppose that $A$ is an $n \times N$ matrix which is depicted from an
SHF$(N;n,m,\{w,r\})$ (The transpose of the matrix).
The elements in $A$ are $1,2,\cdots, m$. Suppose $B$ is the incidence
matrix of a $(w,r;d)$-CFF$(v,m)$. Denote $b_1, b_2, \cdots, b_m$ the rows of
$B$.
We construct a $Nv\times n$ matrix $A'$ by substituting the element $i$ in $A$
by $b_i$. It can be verified that $A'$ is an incidence matrix of a
$(w,r;d)$-CFF$(vN,
n)$. Thus we have the following generalization of \cite{stw}.

\begin{Theorem}
\label{th.uSHF}
If there exists an $(w,r;d)$-CFF$(v,m)$ and an SHF$(N;n,m,\{w,r\})$,
then there exists a
$(w,r;d)$-CFF$(vN,n)$.
\end{Theorem}

Using a recursive method introduced  in \cite{amsw}, the following
result is proved in \cite{stw} (The original result is for
$(w,r)$-CFF. It is straightforward to generalize it to
$(w,r;d)$-CFF.).

\begin{Theorem}
\label{DSrec}
Suppose there exists an $(w,r;d)$-CFF$(N_0,n_0)$,
where $\gcd\left( n_0, (wr)!\, \right) = 1$.
Then there exists an
$(w,r;d)$-CFF$\left(\left( wr +1\right)^k N_0,{n_0}^{2^k}, \right)$
for any integer $k \ge 0$.
\end{Theorem}

To apply the above theorem, we can use the following easy lemma.

\begin{Lemma}
\label{DSdirect}
Suppose that $i$ and $j$ are positive integers
and $n \geq i+j$.
Then there exists an $(i,j)$-CFF$(N,n)$
where $N = \min \left\{ {n \choose i}, {n \choose j} \right\}$.
\end{Lemma}

\p
Let $X$ be an $n$-set.  The set of all
$i$-subsets of $X$ is an $(i,j)$-disjunct system (which is equivalent to an
$(i,j)$-CFF), as is the
set of all $(n-j)$-subsets of $X$.
\qed

Using the above lemma, we can easily obtain the following general result, which
provides an explicitly constructed class of $(w,r)$-CFF
for any $w$ and $r$.

\begin{Theorem}\cite{stw}
\label{DSgenclass}
Suppose that $w$ and $r$ are positive integers.
Let
\[n_0 = \min \{n \geq w+r : \gcd (n, (wr)!\, ) = 1 \}\]
and let $N_0 = \min \left\{ {n_0 \choose w}, {n_0 \choose r} \right\}$.
Then there exists a
$(w,r)$-CFF$\left(  \left( wr +1\right)^k N_0 ,{n_0}^{2^k}\right)$
for any integer $k \geq 0$.
\end{Theorem}

\p
Apply Theorem \ref{DSrec} and Lemma \ref{DSdirect}.
\qed

The following corollary is an immediate consequence of
Theorem \ref{DSgenclass}.

\begin{Corollary}
\label{DSgen}
For any positive integers $w$ and $r$, there exists an
explicit construction for
an infinite class of $(w,r)$-CFF$(N,T)$
in which $T$ is $O\left( (\log N)^{\log (wr + 1)} \right)$.
\end{Corollary}

From  Remark \ref{remark}, we can construct a
$(w,r;d)$-CFF$\left(  \left( wr +1\right)^k (d+1)N_0 ,{n_0}^{2^k}\right)$
for any integer $k \geq 0$, where $N_0$ and $n_0$ are same as in Theorem
\ref{DSgenclass}. Thus we have the following.

\begin{Corollary}
\label{CFgen}
For any positive integers $w$ and $r$, there exists an
explicit construction for
an infinite class of $(w,r;d)$-CFF$(N,T)$
in which $T$ is $O\left( (\log (d+1)N)^{\log (wr + 1)} \right)$.
\end{Corollary}

\subsection{Constructions from codes}

Consider a code $\C$ of length $N$ on an alphabet $Q$
with $|Q| = q$. Then $\C \subseteq Q^N$ and we will call it an
$(N,n,q)$-code if $|\C|= n$.
The elements of $\C$ are called {\it codewords};
each codeword is $x = (x_1, \cdots , x_N)$, where
$x_i \in Q$, $1 \leq i \leq N$.

Suppose $\C$ is an
$(N, n, q)$  code on an alphabet $Q$.
Define $X = \{ 1, \cdots , N\}\times Q$, and for
each codeword $c = (c_1, \cdots , c_N) \in \C$,
define an $N$-subset of $X$ as follows:
\[
B_c = \{(i,c_i) : 1 \le i \le N\}.
\]
Finally, define $\B = \{ B_c : c \in \C\}$. Then we obtain an $N$-uniform
  set system $(X, \B)$. Using this correspondence between a code and a set 
system,
we are able to construct uniform CFF from codes which satisfy certain
properties. The following theorem is easy to prove.

\begin{Theorem}\label{t.CFC}
Suppose that $\C$ is an $(N,n,q)$-code having
minimum distance $D$.
Then there is a $(r;d)$-CFF$(Nq,n)$, where $r=\left\lfloor {N-d-1\over
N-D}\right\rfloor$.
\end{Theorem}

Reed-Solomon codes were used in \cite{eff2,krs,ssw} to construct CFF. Since a
Reed-solomon code is an $(N,q^t,q)$-code with minimum distance $D=N-t+1$, we
have the following.

\begin{Theorem}
\label{RScode}
Suppose $N, q, r$ and $d$ are given, with
$q$ a prime power and $N \le q+1$.   Then there
exists an $(r;d)$-CFF$(qN, T)$ where $T = q^{\lceil (N+r-d-1)/r \rceil}$.
\end{Theorem}

The CFF constructed from a Reed-Solomon code is the same as
 a CFF  constructed
from an orthogonal array.

A shortened Reed-Solomon code is used to construct $r$-CFF in
\cite{dmr}. We can consider this method as a construction of
$(r;d)$-CFF. A shortened Reed-Solomon code is a sub-code of a
Reed-Solomon $(q+1,q^t,q)$-code, which contains all the codewords
whose first $s$ elements are all zeros, where $0 \le s \le t-1$.
So the shortened Reed-Solomon code is a $(q+1-s, q^{t-s},q)$-code.
The minimum  distance of a shortened code is the same as the
original Reed-Solomon code: $D = q-t+2$. So we obtain the
following result.

 \begin{Theorem}
\label{RSshort}
Suppose $q, r$ and $d$ are given, with
$q$ a prime power.   Then there
exists an $(r;d)$-CFF$((q+1-s)q, T)$ where $T = q^{\lceil (q+r-d-s)/r \rceil}$
and $s+d \le q$.
\end{Theorem}

When $r > d+1$, the CFF constructed from a shortened Reed-Solomon code is
better than that
constructed from a Reed-Solomon code. To see that, we may consider the
rate $\log T \over N$.  The rate of a CFF obtained
from a Reed-Solomon code
is
\[
{q+r-d\over rq(q+1)}\log q,
\]
while the rate of the CFF obtained from the shortened code is
\[
 {q+r-d-s\over rq(q+1-s)}\log q.
\]

In \cite{krs}, some algebraic-geometric (Goppa) codes were used to construct 
CFF.
The algebraic-geometric code (or AG code) they used was from \cite{gs}.

Garcia and Stichtenoth \cite{gs} gave an explicit construction of function
fields $F_n$ which are extending from ${\mathbb F}_{q^2}$. They proved that 
$F_n$ has at least
$s(n)+1$ places of degree one and genus at most $g(n)$, where $s(n)=q^{n-1}(q^2
-1)$ and $ g(n)= q^{n-1}(q+1)$. Thus for $n\ge 3$ and $\alpha \le s(n)$ there is
an $(
s(n), q^{2k},q^2)$-code with minimum distance $D \ge s(n) - \alpha$, where $k\ge
\alpha - g(n) +1$.

The following theorem is a generalized form of \cite[Theorem 4]{krs}.

\begin{Theorem}
\label{AGcode}
Suppose $n, q, r$ and $d$ are given, with
$q$ a prime power, $r < q-1$ and $n\ge 3$.   Then there
exists an $(r;d)$-CFF$(q^2s(n), q^{2k})$ where $k \ge \alpha - g(n) +1$ and $r =
{\lceil (s(n)-d-1)/\alpha \rceil}$.
\end{Theorem}

To compare the CFF constructed from Reed-Solomon codes with the CFF constructed
from AG codes, we fix the field ${\mathbb F}_{q^2}, r$ and $d$ for both codes
and consider the rate $\log T\over N$. In the case of Reed-Solomon codes, we let
$N=q^2+1$ and obtain the rate as
\[
{q^2+r-d \over rq^2(q^2+1)}\log q^2.
\]
In the case of Garcia-Stichtenoth codes, the rate is about
\[
{q-r-1\over rq^2(q-1)}\log q^2.
\]

From these rates we can see that if the value of $r$ is much smaller than the
value of $d$, then the CFF obtained from Garcia-Stichtenoth code is better than
the CFF obtained from Reed-Solomon code. Otherwise the CFF obtained from
Reed-Solomon code is more efficient.

\subsection{ Probabilistic methods }

Probabilistic methods were used to investigate the existence of $(w,r)$-CFF
  by many researchers (see \cite{dh,drr,dfft,e,eff2,krs,stw}). All of these
results are similar. Here we record a result provided in
\cite{stw} as follows.

\begin{Theorem}
\label{ds-bound}
Suppose that $T,w$ and $r$ are positive integers.
Define
\[ p = 1 - { w^wr^r \over (w+r)^{w+r} }.\]
Then the following hold:
\begin{enumerate}
\item If
\[ N > { (w+r) \log T \over - \log p},\]
then there exists an $(w,r)$-CFF$(N,T)$.
\item If
\[ N > { (w+r - 1) \log (2T) \over - \log p},\]
then there exists an $(w,r)$-CFF$(N,T)$.
\end{enumerate}
\end{Theorem}

By Remark \ref{remark} we have the following result.

\begin{Theorem}
\label{wrd-bound}
Suppose that $T,w,r$ and $d$ are positive integers.
Define
\[ p = 1 - { w^wr^r \over (w+r)^{w+r} }.\]
Then the following hold:
\begin{enumerate}
\item If
\[ N > { (w+r) \log T \over - (d+1)\log p},\]
then there exists an $(w,r;d)$-CFF$(N,T)$.
\item If
\[ N > { (w+r - 1) \log (2T) \over - (d+1)\log p},\]
then there exists an $(w,r;d)$-CFF$(N,T)$.
\end{enumerate}
\end{Theorem}

Using Chernoff and some probabilistic argument, \cite{sw2} proved the following

theorem about the existence of $k$-uniform generalized
CFF.

\begin{Theorem}
\label{unif-bound}
Suppose that $T,w,r$ and $\ell$ are positive integers.
Define
\[ p = \frac{(\ell - 1)^{r}}{\ell^{w+r-1}}.\]
If 
\[ k > {8\over p}\left( (w+r)\log T -\log w! -\log r!\right),\]
then there exists a $k$-uniform $(w,r;d)$-CFF$(k \ell,T)$ for
\[d = {pk \over 2} + 1.\]
 \end{Theorem}

The following variation is proved in a similar fashion in \cite{sw2}.

\begin{Theorem}
Suppose that $T,w,r$ and $\ell$ are positive integers.
Define
\[ p = \frac{(\ell - 1)^{r}}{\ell^{w+r-1}}.\]
If 
\[ k > {8\over p}\left( (w+r)\log T + t -\log w! -\log r!\right),\]
then the probability that the matrix $A$ is not a
$k$-uniform $(w,r;d)$-CFF$(k \ell,T)$, for
\[d = {pk \over 2} + 1,\]
is at most $e^{-t}$.
\end{Theorem}

\section*{Note added for the submission to arXiv}

This manuscript was written back to 2006. It was accepted by an special issue of Discrete Mathematics.
However, by some reasons, the special issue has never published.  Thanks to Lucia Moura who 
recently suggested
me to submit it to arXiv so that other researchers may access or cite it.

\end{document}